\documentclass[paper=a4, fontsize=12pt, twocolumn]{scrartcl} 

\usepackage[english]{babel} 
\usepackage[protrusion=true,expansion=true]{microtype} 
\usepackage{amsmath,amsfonts,amsthm} 
\usepackage[pdftex]{graphicx} 
\usepackage[svgnames]{xcolor} 
\usepackage[hang, small,labelfont=bf,up,textfont=it,up]{caption}	
\usepackage{subfig}	
\usepackage{fix-cm} 

\usepackage{url}
\usepackage{marginnote}

\usepackage{sectsty} 
\allsectionsfont{
	\usefont{OT1}{phv}{b}{n}
	}

\sectionfont{
	\usefont{OT1}{phv}{b}{n}
	}

\usepackage{fancyhdr} 
    \usepackage{lastpage}	

\fancypagestyle{default}{
    \lhead{}
    \chead{}
    \rhead{}
    \lfoot{\footnotesize Petra Schwer \textbullet ~The meaning of doing mathematics \textbullet ~2025}
    \cfoot{}
    \rfoot{\footnotesize page \thepage\ of \pageref{LastPage}}	

}

\fancypagestyle{last-page}{
    \lhead{}
    \chead{}
    \rhead{}
    \lfoot{\footnotesize Petra Schwer, Professor of Mathematics at Heidelberg University, Germany. \\
    Share your thoughts with me at  \emph{schwer@uni-heidelberg.de} \\
    \url{https://web.mathi.uni-heidelberg.de/ggt/schwer} \\ ~ \\
    The meaning of doing mathematics \textbullet ~2025}
    \cfoot{}
    \rfoot{~ \\~ \\~ \\~ \\
    \footnotesize page \thepage\ of \pageref{LastPage}}	

}

\AtEndDocument{\thispagestyle{last-page}}

\usepackage{lettrine}
\newcommand{\initial}[1]{%
     \lettrine[lines=3,lhang=0.3,nindent=0em]{
     				\color{DarkGoldenrod}
     				{\textsf{#1}}}{}}

\usepackage{titling}

\newcommand{\HorRule}{\color{DarkGoldenrod}	\rule{\linewidth}{1pt}}

\pretitle{\vspace{-30pt} \begin{flushleft} \HorRule 
				\fontsize{40}{40} \usefont{OT1}{phv}{b}{n} \color{Maroon}
                \selectfont 
				}
                
\title{The meaning of \\ doing mathematics \\
\LARGE{{\ }\\Can AI solve all math?}
} 
\posttitle{\par\end{flushleft}\vskip 0.5em}

\preauthor{\begin{flushleft} 
					\large \lineskip 0.5em \usefont{OT1}{phv}{b}{sl} \color{black}}
\author{Petra Schwer, }											
\postauthor{\footnotesize \usefont{OT1}{phv}{m}{sl} \color{Black} 
					Heidelberg, September 2025
					\par\end{flushleft}\HorRule}
\date{}																				

\begin{document}
\maketitle
\thispagestyle{fancy} 			

\initial{'E}\textbf{ventually AI will solve all mathematics.' one of the Deep Mind members said, while we were both sipping at our coffee and nibbling cookies in the common room of Fuld Hall in Princeton earlier this year. \\ 'I don't think that's possible.' I replied and earned a puzzled look. 
}


\medskip

This summer the internet has been full with headlines about AI outperforming Math Olympiad gold medalists or solving Ph.D. level problems in mathematics and hence outsmarting the world's top mathematicians. It is no question, Large Language Models (LLMs for short) are getting more powerful by the day and it is inevitable that society is living through a major phase of change. The world with LLMs is not going to be the same as the one we were used to.  
The elephant in the room is, whether we will eventually be 'replaced by AI'. The question seems, whether AI will solve all mathematics, take away our job, or (even worse so) our purpose. Well, as I said at the coffee break already:  
I don't think that's possible.

In order to see why I arrive at this conclusion, let me discuss  what it is that we actually do, when we do math. 

\subsection*{\textcolor{Maroon}{What do we actually mean by doing math?}}

What drew me to math as a teenager and young adult was its rigor and universal truth.   Math is not debatable. Mathematical statements and reasoning is either wrong or provides the right answer.
Today I am convinced, that if we strip away the humans from mathematics the subject will loose its soul. 
Math is a deeply human endeavor. 

\medskip

We mathematicians spend hours discussing in front of blackboards, sitting at our desks painstakingly checking every detail, choosing just the right notation and names of things in order to formulate and check properties of the objects we look at. 
We end up, hopefully, with a list of statements and their  verification. 
What most of us then produce on a regular basis are papers containing said definitions, theorems, lemmata and proofs thereof. 

But why is it, that we do these things?

\begin{center} \Large
	Doing math is like exploring an unknown landscape.
\end{center}

When standing in front of our blackboards or running computer experiments and simulations we mathematicians are in the middle of an expedition. We set out, often in teams, to explore our mathematical landscapes. Each of us comes with their own experience and partially drawn map in their heads, shaped by what we have been taught, read or heard and  methods we have learned and used. Over time each of us acquires their own individual knowledge and understanding of mathematical objects, their behavior and how they connect. 

Now my map definitely differs from yours. So already when standing in front of those blackboards with just one other colleague we often struggle to find common words for the patterns we observe and the structures we want to describe. 

\begin{center} \Large
	My mental map \\ of mathematics\\ is useless in case I can not communicate it well with others. 
\end{center}

In fact, a large part of the work we do is to curate and illustrate pieces of our own mental map. This includes developing methods, defining properties and objects, identifying statements, writing proofs and much more. 
We spend all those hours trying to develop the best possible notions and language in order to communicate clearly what we see and understand with other mathematicians or with  scientists in other fields.

\smallskip

\subsection*{\textcolor{Maroon}{How do we communicate mathematics among mathematicians?}}

Over the course of a long history mathematicians have acquired powerful tools that allow us to communicate mathematics with others. 

We use symbols\footnote{Did you know that the equal sign '=' only became popular in the early 17th century? Several competing signs were used to prior to it. Among them the symbol nowadays used to indicate parallelism of lines.} or words, like \emph{natural numbers}, to abbreviate concepts lengthy to describe. 
We try to use suggestive names for objects and (maybe even more importantly) try to find well chosen notation that facilitates  proofs and computations. We structure written proofs in an easy to follow way. We introduce dictionaries to translate a statement from one point of view (say a geometric one) to a different viewpoint (an algebraic one, for example). 

\begin{center} \Large
	Learning math \\ also involves learning about \\
	communication tools\\ and \\ how to use them properly. 
\end{center}

The most powerful among those tools is \emph{proofs}. We deduce, with logical arguments, one mathematical statement from other, known statements. 

Proofs started with Euclid about 2300 years ago and had their foundational crisis 
culminating in Gödel's famous incompleteness theorems 
at the beginning of the 20th century. Bourbaki has then shaped how much symbols and rigor we expect and most recently proof assistants influence how formally checkable a statement may be. 

What is considered a proof
is based on shared communal knowledge and a common understanding of formalism and rigor.
It will depend on the audience how detailed a proof needs to be or which tools one will use.
Proofs are the backbone of modern mathematics, yet what is considered a valid proof 
has changed over time and will likely change in the future. 


Having a proof (or claiming to have one) is quite useless to the community, though, unless we can communicate our understanding to others. We don't want to believe statements, we want to understand how and why the statement connects to the map we have in our minds. 

\begin{center} \Large
It is a bit like in ballet: \\ You put in a lot of work to make it look easy in the end. 
\end{center}

I once received a referee report on an expository piece of work saying "I like what the author did to Theorem X - it looks almost natural now." I think this is the whole point of writing: To make the content look natural. Much like ballet, which it is not just about performing a given sequence of complicated moves but also about making it look beautiful, light and easy.

\subsection*{\textcolor{Maroon}{What is it all about then?}}

We long for precision and seek eternal truth. And we want to ensure we did not make mistakes on our way. Mathematicians as a whole have come a long way in establishing a language and tools that allows us to actually sometimes succeed in doing so. However, true statements that nobody understands are useless. 

\begin{center} \Large
   Truth needs meaning. 
\end{center}

We want to make it possible for others to see pieces of the mathematical world as we do and to use it for their own undertakings. Communicating mathematics in a way that empowers others is an interpersonal matter.

\smallskip

Mathematicians prove theorems and write papers. One may thus think that the most valuable contribution of a mathematician is a list of proven-to-be-true statements. As a community we reinforce this viewpoint by how we measure success, evaluate grants and hire people. 

And yet, not everything is visible in those papers. In fact, mathematicians have developed the habit of not writing about those things we investigated (during the search for a proof) that remained irrelevant for the final argument. But it may be exactly these seemingly useless explorations that lay the foundation for future work, new questions and new mathematics to emerge.  

\smallskip

Writing a good piece of mathematics is a creative process sharing similarities with the making of a piece of art or music. We provide \emph{elegant} arguments and admire \emph{beautiful} proofs. And sometimes (not everybody shares the same taste) it is true joy to read or listen to the ideas and findings of another mathematician - maybe a bit like listening to a well crafted piece of music.

As Bill Thurston\footnote{He wrote this in his reply to the  mathoverflow question titled \emph{What's a mathematician to do?}. I strongly recommend reading his entire answer here: \url{https://mathoverflow.net/q/44213}.} put it: \emph{The product of mathematics is clarity and understanding. }

\medskip

Reducing math to solving problems misses the point of what we\footnote{Obviously this article is not written by all mathematicians. In the past few months I have discussed various aspects of the content of this writing with many colleagues. Several of them expressed similar thoughts and agreed with me on the necessity to communicate them. I am hence confident in speaking of and with the voice of a larger 'we' here. } do. And it misses the point of what makes mathematics attractive and worthwhile doing for me (and many of my colleagues) as well as for society.

\smallskip

My master thesis advisor often claimed that one is not trained but raised as a mathematician. When working in industry for a bit with a very diverse crowd of colleagues, I got an idea what this statement means. Mathematicians think differently. We have a very different way of approaching and structuring hard problems than people with other backgrounds and training. Not because math is superior, but because we have been trained, for years, to think with the most open mind possible and to approach problems fearlessly. 

\begin{center} \Large
	Math is a way of thinking. 
\end{center}

{This is one of the reasons why other professions and a diverse range of industries benefit from talking to or working with mathematicians despite the mathematician's initial lack of training on the actual topic.}

\subsection*{\textcolor{Maroon}{Now, can AI solve all mathematics?}}

    I don't think so. 
    First of all, mathematics is not finite. 
    There is simply no such thing as 'all math'.  
    More importantly: 
    
    \begin{center} \Large
        Math is not just \\ a list of problems\\ to be solved. 
    \end{center}

    {\tiny We've done a great job in teaching generations of highschool students quite the opposite of this, unfortunately.} \\

	Mathematicians design and map landscapes, solve and set problems, ask questions and aim to find structure in what we explore. And we refine and use our tools to share our discoveries with others. 
    
    \smallskip
    
    'Solving problems' is a task the existing LLMs excel in. By design, a large language model is good at repeating patterns available in their training data. If you ask a question it will predict what the most likely answer to this string of words will be.  
    This is why, for example, LLMs are so good at solving math olympiad puzzles. There is a lot of training data freely available on the internet.        
    LLMs can also be used to help find connections between seemingly unrelated areas. Machines are much quicker in connecting the dots over long distances and thus discover connections that exist in writing but which humans would need a (very) long time to see. 

    AI has succeeded and will continue to succeed in outperforming mathematicians in solving problems. 
    Some mathematicians are involved in challenging the latest LLMs. Collections of open questions and problems, like FrontierMath\footnote{See \url{https://arxiv.org/abs/2411.04872}.}, are being designed to test the boundaries of the current models.   
	
    It is not the question whether AI will (help) solve mathematical problems. 

    \begin{center} \Large
        AI does and will shape\\ how we work. 
    \end{center}
    
    We are already seeing and living these changes. AI might provide powerful tools that, especially when combined with formal proof assistants, might continue to change the way true statements can be produced. 
    In a best case scenario AI will be our travel companion. 
	
    It is  not for the first time society sees radical changes like this. Mankind has invented the wheel, books, the internet, ... 
    While some professions are drastically challenged for the first time, mathematics has a head start on dealing with machines that are faster than us in performing algorithms and formal computations. This goes from abaci, slide rulers, calculators to turing machines and proof assistants. 
   
    \begin{center} \Large
        The question is\\ how we want to use AI. 
    \end{center}
    
    The current developments may suggest that we can jump over the hard parts and just curate the output AI gives us. This would be much like  running a marathon while skipping all of the shorter training sessions. It could well be that the way we train\footnote{We train ourselves and our students. The effect AI will have on our teaching is not discussed here. } is going to change.  
    We still teach the multiplication tables to first graders, though. Despite the existence of all the aforementioned devices. 

	\vfill\eject 

    To me the question is also: how can we make sure that with the changes we see the human part of doing mathematics is not getting lost along the way?

\subsection*{\textcolor{Maroon}{What I want you to remember}}

One thing I learned over coffee this spring is, that some people seem to think mathematics is just about solving problems. 

If we do not want to be reduced to this  
we need to {remember} that the point of  

    \begin{center} \Large
        'doing mathematics'\\ is to explore, build and communicate\\ a universal, everlasting and  true\\ map of mathematics\\  that other humans understand and are capable of\\ putting into use. 
    \end{center}

And, maybe even more importantly, we need to let people know, that this is what we do.

\vspace{30ex}

\end{document}